\documentclass[12pt]{amsart}
\usepackage{amsbsy}
\usepackage{graphicx,epsfig,subfigure,psfrag}
\begin{document}

\newtheorem{theorem}{Theorem}
\newtheorem{proposition}{Proposition}
\newtheorem{lemma}{Lemma}
\newtheorem{corollary}{Corollary}
\newtheorem{definition}{Definition}
\newtheorem{remark}{Remark}
\newcommand{\tex}{\textstyle}
\numberwithin{equation}{section} \numberwithin{theorem}{section}
\numberwithin{proposition}{section} \numberwithin{lemma}{section}
\numberwithin{corollary}{section}
\numberwithin{definition}{section} \numberwithin{remark}{section}
\newcommand{\ren}{\mathbb{R}^N}
\newcommand{\re}{\mathbb{R}}
\newcommand{\n}{\nabla}
\newcommand{\iy}{\infty}
\newcommand{\pa}{\partial}
\newcommand{\fp}{\noindent}
\newcommand{\ms}{\medskip\vskip-.1cm}
\newcommand{\mpb}{\medskip}
\newcommand{\AAA}{{\bf A}}
\newcommand{\BB}{{\bf B}}
\newcommand{\CC}{{\bf C}}
\newcommand{\DD}{{\bf D}}
\newcommand{\EE}{{\bf E}}
\newcommand{\FF}{{\bf F}}
\newcommand{\GG}{{\bf G}}
\newcommand{\oo}{{\mathbf \omega}}
\newcommand{\Am}{{\bf A}_{2m}}
\newcommand{\CCC}{{\mathbf  C}}
\newcommand{\II}{{\mathrm{Im}}\,}
\newcommand{\RR}{{\mathrm{Re}}\,}
\newcommand{\eee}{{\mathrm  e}}
\newcommand{\LL}{L^2_\rho(\ren)}
\newcommand{\LLL}{L^2_{\rho^*}(\ren)}
\renewcommand{\a}{\alpha}
\renewcommand{\b}{\beta}
\newcommand{\g}{\gamma}
\newcommand{\G}{\Gamma}
\renewcommand{\d}{\delta}
\newcommand{\D}{\Delta}
\newcommand{\e}{\varepsilon}
\newcommand{\var}{\varphi}
\renewcommand{\l}{\lambda}
\renewcommand{\o}{\omega}
\renewcommand{\O}{\Omega}
\newcommand{\s}{\sigma}
\renewcommand{\t}{\tau}
\renewcommand{\th}{\theta}
\newcommand{\z}{\zeta}
\newcommand{\wx}{\widetilde x}
\newcommand{\wt}{\widetilde t}
\newcommand{\noi}{\noindent}
\newcommand{\uu}{{\bf u}}
\newcommand{\xx}{{\bf x}}
\newcommand{\yy}{{\bf y}}
\newcommand{\zz}{{\bf z}}
\newcommand{\aaa}{{\bf a}}
\newcommand{\cc}{{\bf c}}
\newcommand{\jj}{{\bf j}}
\newcommand{\ggg}{{\bf g}}
\newcommand{\UU}{{\bf U}}
\newcommand{\YY}{{\bf Y}}
\newcommand{\HH}{{\bf H}}
\newcommand{\GGG}{{\bf G}}
\newcommand{\VV}{{\bf V}}
\newcommand{\ww}{{\bf w}}
\newcommand{\vv}{{\bf v}}
\newcommand{\hh}{{\bf h}}
\newcommand{\di}{{\rm div}\,}
\newcommand{\ii}{{\rm i}\,}
\newcommand{\inA}{\quad \mbox{in} \quad \ren \times \re_+}
\newcommand{\inB}{\quad \mbox{in} \quad}
\newcommand{\inC}{\quad \mbox{in} \quad \re \times \re_+}
\newcommand{\inD}{\quad \mbox{in} \quad \re}
\newcommand{\forA}{\quad \mbox{for} \quad}
\newcommand{\whereA}{,\quad \mbox{where} \quad}
\newcommand{\asA}{\quad \mbox{as} \quad}
\newcommand{\andA}{\quad \mbox{and} \quad}
\newcommand{\withA}{,\quad \mbox{with} \quad}
\newcommand{\orA}{,\quad \mbox{or} \quad}
\newcommand{\ef}{\eqref}
\newcommand{\ssk}{\smallskip}
\newcommand{\LongA}{\quad \Longrightarrow \quad}
\def\com#1{\fbox{\parbox{6in}{\texttt{#1}}}}
\def\N{{\mathbb N}}
\def\A{{\cal A}}
\newcommand{\de}{\,d}
\newcommand{\eps}{\varepsilon}
\newcommand{\be}{\begin{equation}}
\newcommand{\ee}{\end{equation}}
\newcommand{\spt}{{\mbox spt}}
\newcommand{\ind}{{\mbox ind}}
\newcommand{\supp}{{\mbox supp}}
\newcommand{\dip}{\displaystyle}
\newcommand{\prt}{\partial}
\renewcommand{\theequation}{\thesection.\arabic{equation}}
\renewcommand{\baselinestretch}{1.1}
\newcommand{\Dm}{(-\D)^m}

\title
{\bf Formation of shocks in \\  higher-order nonlinear dispersion
PDEs:\\
  nonuniqueness
 and nonexistence of  entropy}

\author {
Victor A.~Galaktionov}

\address{Department of Mathematical Sciences, University of Bath,
 Bath BA2 7AY, UK}
\email{vag@maths.bath.ac.uk}


  \keywords{Odd-order quasilinear PDE, shock waves,
  self-similar solutions, blow-up, nonuniqueness.
  {\bf Submitted to:} J. Differ. Equat.
  }
 \date{\today}

\begin{abstract}

 Formation of  shocks  for the third-order nonlinear
dispersion PDE
 \be
 \label{01}
 u_t=(uu_x)_{xx} 
 \quad \mbox{in}
 \quad \re \times \re_+,
 \ee
 is studied by construction various self-similar solutions exhibiting gradient blow-up in finite time,
 as $t \to T^- \in (0,\iy)$, with locally bounded final time profiles $u(x,T^-)$.
 These are shown to   admit infinitely many discontinuous shock-type similarity
 extensions for $t>T$, all of them satisfying
 generalized Rankine--Hugoniot's condition at shocks. As a result,  the principal
 nonuniqueness of solutions  of the Cauchy problem after blow-up time is inherited, and,
 under certain hypothesis,
  any sufficiently general ``entropy-like" approach to such
 problems to detect a unique continuation beyond singularity becomes
 illusive.
 In other words,  any attempt to create ``entropy theory" for
 (\ref{01}) along the lines of the classic one for
   scalar conservation laws
  (developed by Oleinik and Kruzhkov
in the 1950-60s), such as Euler's equation
 $u_t + uu_x=0$
 in $\re \times \re_+$,
 or by using other more recent
  successful  ideas for hyperbolic systems, is hopeless,
  regardless the fact that
  these PDEs have  a number of  similar
  features of formation of
   shock and rarefaction
  waves.

\end{abstract}

\maketitle






\section {Introduction: nonlinear dispersion PDEs and main directions of study}
 \label{Sect1}

We study  ``micro-structure" of shock-type finite time
singularities that can occur in higher-order nonlinear dispersion
equations (NDEs). To describe key features of formation of such
single point singularities, it suffices
 to consider the NDE of the minimal third order.
 Its study
eventually leads to quite a pessimistic conclusion concerning
uniqueness of ``entropy" solutions and on nonexistence of any
hypothetical entropy.


 \subsection{NDEs: nonlinear dispersion equations}

 In the present paper,
 we continue our study  began in \cite{GalEng}--\cite{GPnde} of
   basic aspects of singularity formation and  approaches
   to
existence-uniqueness-entropy  for odd-order {\em nonlinear
dispersion} (or {\em dispersive}) equations.
The simplest  canonical  model 
is the 
{\em third-order quadratic NDE} (the NDE--$3$)
 \be
 \label{1}
  \mbox{$
 u_t=
 (uu_x)_{xx}  \equiv u u_{xxx}+ 3 u_x u_{xx}
 \quad \mbox{in}
 \quad \re \times (0,T), \,\,\, T >0.
  $}
 \ee
We pose for (\ref{1}) the Cauchy problem (the CP)
  with locally bounded and
integrable initial data
 \be
 \label{2}
 u(x,0)=u_0(x) \quad \mbox{in} \quad \re.
  \ee
 It is principal that we consider the CP, where the solution $u(x,t)$ is supposed to
 be defined by initial data \ef{2} only,
 plus, of course a generalized Rankine--Hugoniot-type condition on the speed
 of propagation of shocks, which follows from the equation
 integrated in a shock neighbourhood.
 In other words, we are assuming that the CP
 does not require any
  {\em a priori} posed conditions on the shock wave lines (though
  of course such ones exist and can be determined {\em a  posteriori}).
 Otherwise, with such conditions, we arrive at a {\em free
 boundary problem} (an FBP) for the NDE \ef{1}, which requires
 other mathematical methods of study and can be well-posed (unlike
 the CP). We will  touch possible FBP settings for
 \ef{1}.

The physical motivation of the NDEs such as \ef{1} and other
odd-order nonlinear PDEs, which appear in many areas of
application, with a large number of key references,  are available
in surveys in \cite[\S~1]{GalNDE5} or \cite[\S~1]{GPnde}, so we
omit any discussion on these applied issues.


\subsection{NDEs in general PDE theory}

For future convenience, we briefly present a more mathematical
discussion around NDEs; cf. \cite[\S~1]{GPnde}. In the framework
of general theory of nonlinear evolution PDEs of the first order
in time, the NDE (\ref{1}) appears the third among
other  canonical evolution quasilinear degenerate equations:
 \be
 \label{ca1}
  \mbox{$
 u_t= - \frac 12\, (u^2)_x \quad (\mbox{the conservation law}),
 \qquad\qquad\qquad\qquad\qquad\,\,\,\,
  $}
  \ee
 \be
 \label{ca2}
  \mbox{$
 u_t= \frac 12 \,(u^2)_{xx}\quad (\mbox{the porous medium equation}),
 \qquad\qquad\qquad\quad\,\,\,\,
 $}
 \ee
 \be
\label{ca3} \mbox{$
 u_t= \frac 12 \,(u^2)_{xxx}\quad (\mbox{the nonlinear dispersion equation}),
 \qquad\qquad\,\,\,\,\,\,\,
 $}
 \ee
  \be
\label{ca4}
 \mbox{$
 u_t= -\frac 12 \,(|u|u)_{xxxx}\quad (\mbox{the $4^{\rm th}$-order nonlinear diffusion
 equation}).
 $}
 \ee
 In (\ref{ca4}), the quadratic nonlinearity $u^2$ is replaced by
 the monotone one $|u|u$ in order to keep the parabolicity on
 solutions of changing sign. The same can be done in the PME
 (\ref{ca2}), though this classic version is parabolic on
 nonnegative solutions, a property that is preserved by the
 Maximum Principle.
Further extensions of the list by including
  \be
   \label{ca5}
    \mbox{$
 u_t= -\frac 12 \,(u^2)_{xxxxx}\quad (\mbox{the NDE--5}) \quad \mbox{and}
 \qquad\qquad\qquad\qquad\qquad\,\,\,\,\quad \quad
 $}
 \ee
 \be
\label{ca6}
 \mbox{$
 u_t= \frac 12 \,(|u|u)_{xxxxxx}\quad (\mbox{the $6^{\rm th}$-order nonlinear diffusion
 equation), etc.,}
 $}
 \ee
are not that essential. These PDEs belong to the same families as
(\ref{ca3}) and (\ref{ca4}) respectively with similar covering
mathematical concepts (but more difficult in some details).

 Mathematical theory of the first two equations, (\ref{ca1}) (see
detailed survey and references below) and (\ref{ca2}) (for quoting
 main achievements of  PME theory developed in the 1950--80s,
see e.g., \cite[Ch.~2]{AMGV}), was essentially completed in the
twentieth century. It is curious that looking more difficult the
fourth-order nonlinear diffusion equation (\ref{ca4}) has a
monotone operator in $H^{-2}$, so the Cauchy problem admits a
unique weak solution as follows from classic theory of monotone
operators; see Lions \cite[Ch.~2]{LIO}. Of course, some other
qualitative properties of solutions of (\ref{ca4}) are  more
difficult and remain open still.

It turns out that, rather surprisingly,  the third order NDE
(\ref{ca3}) is the only one in this basic list that has rather
obscure understanding and lacking of a reliable mathematical basis
concerning generic singularities, shocks, rarefaction waves, and
entropy-like theory.

\subsection{Mathematical preliminaries: analogies with conservation laws, Riemann's problems,
 and earlier results}
 \label{Sect1.2}

  As a key feature of our analysis,
 equation
(\ref{1}) inherits clear similarities 
of the behaviour for the first-order conservation laws such as 1D
{\em Euler's equation} (same as (\ref{ca1})) from gas dynamics
 \be
 \label{3}
  u_t + uu_x=0
\quad \mbox{in} \quad \re\times \re_+,
 \ee
 whose entropy theory  was created by Oleinik \cite{Ol1, Ol59} and Kruzhkov
\cite{Kru2} (equations in $\ren$) in the 1950--60s; see details on
the history, main results, and modern developments in the
well-known monographs \cite{Bres, Daf, Sm}\footnote{First study of
discontinuous  shocks for  quasilinear equations
  was performed  by Riemann in 1858 \cite{Ri58} (by Riemann's method); see
\cite{Chr07, Pom08} for details. The implicit solution of the
problem \ef{3}, $u=u_0(x-u \, t)$ (containing the key
 wave ``overturning" effect), was obtained earlier by Poisson in 1808
\cite{Poi08}; see \cite{Pom08}.}.

\ssk

As for (\ref{3}), in view of the full divergence of the equation
(\ref{1}), it is natural to define weak solutions. For
convenience, we present here  a standard definition, mentioning
however that, in fact, the concept of weak solutions for NDEs even
in fully divergent form is not entirely consistent and/or very
helpful,  to say nothing about other non-divergent equations
admitting no standard weak formulation at all; see \cite{GPnde}.


\ssk

\noi{\bf Definition \ref{Sect1}.1.} {\em A function $u=u(x,t)$ is
a
weak solution of $(\ref{1})$, $(\ref{2})$ if} 

 (i) $u^2 \in L^1_{\rm loc}(\re\times (0,T))$,

 (ii) {\em $u$ satisfies $(\ref{1})$
in the weak sense:
for any test function $\varphi(x,t) \in C_0^\infty(\re \times
(0,T))$},
 \be
 \label{21}
 \mbox{$
 \iint u \varphi_t=\frac 12\,  \iint u^2 \varphi_{xxx},
 $}
  \ee
{\em and} (iii) {\em satisfies the initial condition $(\ref{2})$
in the sense of distributions},
 \be
 \label{22}
 \mbox{$
{\rm ess} \,  \lim_{t \to 0}\int u(x,t) \psi(x) = \int u_0(x) \psi
(x) \quad
 \mbox{\em for any} \quad \psi \in C_0^\infty(\re).
  $}
   \ee

The assumption $T<\infty$ is often essential, since, unlike
(\ref{3}), the NDE (\ref{1}) can produce complete blow-up from
bounded data, \cite[\S~4]{GPnde}.  Thus, again similar to
(\ref{3}), one observes a typical difficulty: according to
Definition \ref{Sect1}.1,
 both discontinuous step-like functions
  \be
  \label{Ri1}
  S_\mp(x) = \mp{\rm sign}\, x = \mp \left\{
 \begin{matrix}\,\,\,
 1 \,\,\, \mbox{for} \,\,\, x>0,\\
 -1
 \,\, \mbox{for}
\,\,\, x < 0,
 \end{matrix}
 \right.
  \ee
 are weak {\em stationary}  solutions of (\ref{1}) satisfying
 \be
 \label{Ri33}
 (u^2)_{xxx}=0
  \ee
  in the  weak sense,
 since in (\ref{21}) $u^2(x)=S_\mp^2(x) \equiv 1$ is $C^3$
smooth ($x=0$ does not count). Again referring to entropy theory
for conservation laws (\ref{3}),  it is well-known that
 \be
 \label{rr1}
 \begin{matrix}
u_-(x,t) \equiv S_-(x) \,\,\, \mbox{is an entropy shock wave,
and}\ssk\\
 u_+(x,t) \equiv S_+(x) \,\,\, \mbox{is not an entropy solution}.\qquad
 \end{matrix}
  \ee
This means that
 \be
 \label{rr41}
 u_-(x,t) \equiv S_-(x)=-{\rm sign} \, x
  \ee
is the unique entropy solution of the PDE (\ref{3}) with the same
initial data $S_-(x)$. On the contrary, taking $S_+$ initial data
  yields the   {\em rarefaction wave}
    with a simple similarity piece-wise linear
    structure
   \be
   \label{rr6}
   \mbox{$
    u_0(x)=S_+(x)={\rm sign} \, x \,\,\, \Longrightarrow \,\,\,
   u_{+}(x,t)= F(\frac xt) = \left\{ \begin{matrix}
    -1 \quad \mbox{for} \,\,\, x<-t, \cr
 \,\,\frac xt \quad  \mbox{for} \,\,\, |x|<t, \cr
  1 \quad \,  \mbox{for} \,\,\, x>t.
 \end{matrix}
   \right.
    $}
 \ee

The questions on  justifying the same classification of main two
Riemann's problems with data (\ref{Ri1}) for the NDE (\ref{1}) and
to construct the corresponding rarefaction wave for $S_+(x)$, as
an analogy of (\ref{rr6}) for the conservation law were addressed
in \cite{GPnde, GPndeII}. This was done by studying the following
self-similar solutions of \ef{1}:
 \be
 \label{2.1}
  \tex{
 u_-(x,t)=f(y), \quad y= \frac x{(-t)^{1/3}}
  \LongA
 (f f')''= \frac 13 \, f'y
 \quad \mbox{in} \quad \re, \quad f(\mp
  \infty)=\pm 1.
  }
  \ee
   Here, by translation, the blow-up time in (\ref{2.1}) reduces to
   $T=0$.
   It was shown that,
 in the sense of distributions or in
 $L^1_{\rm loc}$,
 \be
 \label{2.3}
  u_-(x,t) \to S_-(x) \quad \mbox{as}
  \quad t \to 0^-.
   \ee
 Therefore, $S_-(x)$ is a $\d$-entropy shock wave (see  concepts
 in
 \cite{GPndeII}), while  $S_+(x)$ is not and creates a typical {\em
 rarefaction wave} given by the global similarity solution
  \be
 \label{2.1NN}
  \tex{
 u_+(x,t)=F(y), \quad y= \frac x{t^{1/3}}
  \LongA
 (F F')''=- \frac 13 \, F'y
 \quad \mbox{in} \quad \re, \quad F(\mp
  \infty)=\mp 1,
  }
  \ee
  where indeed $F(y) \equiv -f(y)$, so that, in the same sense as
  in \ef{2.3},
   $$
u_+(x,t) \to S_+(x) \quad \mbox{as}
  \quad t \to 0^+.
   $$

\subsection{Layout of the paper: new shock  patterns and
nonuniqueness}

In Section \ref{S2}, we study new shock patterns, which are
induced by other similarity solutions:
 \be
 \label{s1}
  \tex{
  u_{-}(x,t)=(-t)^\a f(y), \quad y = \frac x{(-t)^\b}, \quad
  \b= \frac {1+\a}3, \quad \mbox{where $\a \in (0, \frac 12)$ and}
   }
   \ee
   \be
   \label{s2}
    \tex{
    (ff')''- \b f'y + \a f=0 \inB \re_-, \quad f(0)=f''(0)=0.
     }
     \ee
The anti-symmetry conditions in \ef{s2} allow to extend the
solution to the positive semi-axis  $\{y>0\}$ by $-f(-y)$ to get a
global pattern. The case $\a<0$ in \ef{s1}, corresponding to the
strong complete blow-up was studied in detail in
\cite[\S~4]{GPnde} in the parameter range
 \be
 \label{ra1}
  \tex{
  \a \in \big[-\frac 1{10},0\big) \whereA a_{\rm c}= -\frac 1{10}
  \,\,\, \mbox{is a special critical exponent}.
  }
   \ee
   For convenience,
we  revise some of these blow-up results in the range \ef{ra1}
 in Section \ref{S2}.

Obviously, the solutions \ef{2.1}, which are suitable for
Riemann's problems, correspond to the simple case $\a=0$ in
\ef{s1}. We prove that, using positive $\a$, allow to get first
{\em gradient blow-up} at $x=0$ as $t \to 0^-$, as a weak
discontinuity, where the final time profile remains locally
bounded and continuous:
 \be
 \label{s3}
  u_{-}(x,0^-)= \left\{
   \begin{matrix}
   C_0 |x|^{\frac \a \b}\,\,\,\, \forA x<0, \ssk\\
- C_0 |x|^{\frac \a \b} \forA x>0,
 \end{matrix}
 \right.
  \ee
  where $C_0>0$ is an arbitrary constant.
  Note that
 $\frac \a \b<1$ for $\a < \frac 12$.

   Therefore, the wave braking
  (``overturning") begins at $t=0$, and,  in Section \ref{S3}, we show that it is
  performed again in a self-similar manner and is described by
  similarity solutions
\be
 \label{s1N}
  \tex{
  u_{+}(x,t)=t^\a F(y), \quad y = \frac x{t^\b}, \quad
  \b= \frac {1+\a}3, \quad \mbox{where}
   }
   \ee
   \be
   \label{s2N}
    \left\{
    \begin{matrix}
    (F F')''+ \b F'y - \a F=0 \inB
    \re_-,\qquad\qquad\qquad\qquad\quad\,\,
  \ssk\\
    F(0)=F_0>0, \,\,\,   F(y) = C_0
    |y|^{\frac \a \b}(1+o(1)) \asA y \to - \iy,
     \end{matrix}
     \right.
     \ee
     where the constant $C_0>0$ is fixed by blow-up data
     \ef{s3}.
The asymptotic behaviour as $y \to -\iy$ in \ef{s2N} guarantees
the continuity of the global discontinuous pattern (with $F(-y)
\equiv -F(y)$) at the singularity blow-up instant $t=0$, so that
 \be
 \label{s3N}
 u_-(x,0^-)= u_+(x,0^+) \inB \re.
  \ee

Rather surprisingly, for a fixed $C_0>0$
 in \ef{s3} obtained by blow-up evolution as $t \to 0^-$, we
 find infinitely many solutions of the extending problem
 \ef{s2N}. This
 family of solutions $\Phi_{C_0}=\{F(y;F_0), \,\, F_0>0\}$
is a one-dimensional curve parameterized by arbitrary
 \be
 \label{s4N}
 F_0=F(0)>0,
  \ee
which does not have any clear boundary or extremal points.
 We also show that for $F_0=0$ \ef{s2N} does not have a solution.
 In other
words, this solution family does not contain any ``minimal",
``maximal", or ``extremal" points in any reasonable  sense, which
might
play a role of a unique ``entropy" one chosen by
introducing a hypothetical entropy inequalities, conditions, or
otherwise.

A first immediate consequence of our similarity blow-up/extension
analysis is as follows:
 \be
 \label{s5N}
 \fbox{$
 \mbox{in the CP, formation of shocks for the NDE (\ref{1}) leads to
 nonuniqueness.}
 $}
 \ee
The second conclusion is more subtle and is based on the mentioned
above fact on the homogeneous structure of the functional set
$\Phi_{C_0}$: if $\Psi_{C_0}=\{u_+(x,t), \, F \in \Phi_{C_0}\}$ is
the whole set of weak solutions of \ef{1} with initial data
\ef{s3}, then, for the Cauchy problem for \ef{1},
 \be
 \label{s6N}
  \fbox{$
 \mbox{there exists
  no  general ``entropy mechanism" to choose a unique
 solution.}
  $}
  \ee
Of course, we cannot exclude a hypothetical situation, when there
exists another, non-similarity solution of this problem, which
thus does not belong to the family $\Psi_{C_0}=\{u_+\}$ and is the
right solution with a proper entropy-type specification. In our
opinion, this is suspicious, so that we claim that \ef{s5N} and
\ef{s6N} show that the problem of uniqueness of weak solutions for
the NDEs such as \ef{1} cannot be solved in principal. On the
other hand, in a FBP setting by adding an extra suitable
 condition on shock lines, the problem might be well-posed with a
 unique solution, though proofs can be very difficult.

In other words, {\em non-uniqueness} in the CP is a non-removable
issue of PDE theory for higher-order degenerate nonlinear
odd-order equations (and possibly not only those).  In fact, the
non-uniqueness of solutions of \ef{s2N} has a pure and elementary
dimensional nature. Indeed, this is a third-order ODE with the
general solution depending on three parameters, which are more
than and excessively  enough to shoot the right behaviour at
infinity and hence to get the continuity \ef{s3N} at the blow-up
time. This guarantees non-uniqueness of the extension for $t>0$
and is a core of this difficulty.

Increasing the order of the PDE under consideration, we then
enlarge  the dimension of the parametric space, and this surely
will imply even stronger non-uniqueness conclusions. For example,
non-uniqueness and non-entropy features are available for the {\em
fifth-order nonlinear dispersion equation} (the NDE--5)
 \be
 \label{nd1}
   u_t=- (u u_x)_{xxxx} \inB \re \times \re_+
   \ee
   and  others;
   see basic models in
   \cite{GPndeII, GalNDE5}.

\section{Gradient blow-up similarity solutions}
 \label{S2}

 In this section, we consider the blow-up ODE problem \ef{s2}.
Actually, this ODE is not that difficult for application of
standard shooting methods, which  in greater detail are explained
in \cite[\S~3,4]{GPnde}. Moreover, even similar fifth-order ODEs
associated with the shocks for the NDE--5 \ef{nd1} also admit
similar shooting analysis \cite{GalNDE5}, though, in view of the
essential growth of the dimension of the phase space,  some more
delicate issues on, say, uniqueness of certain orbits, become very
difficult or even remain open. Therefore, in what follows, we will
omit or even will not mention some technical details concerning
the problem \ef{s2} and recommend consulting \cite{GPnde} in case
of troubles. We will widely use numerical methods for illustrating
and even justifying some of our conclusions. For the third-order
equations such as (\ref{1}), this and further numerical
constructions are performed by the {\tt MatLab} with the standard
{\tt ode45} solver therein; see more details in
\cite[\S~3,4]{GPnde}.

Let us begin with a simpler fact concerning this problem and some
simple asymptotics for matching purposes. We recall the elementary
symmetry of the ODE \ef{s2}
 \be
     \label{symm88}
      \left\{
      \begin{matrix}
     f \mapsto -f,\\
     y \mapsto -y,
      \end{matrix}
      \right.
      \ee
      which allows us to put two conditions at the origin.
      Such solutions have
 a sufficiently  regular asymptotic expansion near
the origin: for any $A <0$, there exists a unique solution of the
ODE (\ref{s2}), satisfying
 \be
 \label{2.7}
  \mbox{$
 f(y) = A y+ \frac {1-2\a}{72}\, y^3+
  \frac {(1-2\a)^2}{72^2}\, \frac 1A \, y^5 +...\, \asA y \to 0.
  $}
  \ee
The uniqueness of such  asymptotics is traced out by using
Banach's  Contraction Principle
 applied to the equivalent integral equation in the metric
 of $C(-\d,\d)$, with $\d>0$  small.

 We  use the following scaling invariance of the
ODE in (\ref{s2}): if $f_1(y)$ is a solution, then
 \be
 \label{2.8}
 \mbox{$
 f_a(y) = a^3 f_1\big(\frac y a\big) \quad \mbox{is a solution for any $a \not =
 0$}.
  $}
  \ee
By scaling \ef{2.8}, the parameter $A<0$ can be reduced to a
single value, say $A=-1$.

Let us describe the necessary bundles of the 3D (in fact, 4D, for
the non-autonomous case) dynamical system \ef{s2}. Firstly, due to
the scaling symmetry \ef{2.8}, there exists the explicit solution
 \be
 \label{2.9}
  \mbox{$
  f_{*}(y)= \frac 1{60} \, y^3<0 \quad \mbox{for} \quad y < -1.
 $}
   \ee
  The overall bundle about \ef{2.9} is 2D, which is obtained by
   the linearization:
   \be
   \label{kk100}
   f(y)=f_*(y)+ Y(y) \LongA  [f_*(y) y]'''- \b Y' y+ \a Y=...=0.
    \ee
The linearly independent solutions are
\be
\label{kk2}
 Y(y)= y^m \LongA (m-3)(m^2+9m -20 \a -2)=0.
  \ee
Hence there exist two negative roots $m _\pm$ composing a 2D
stable manifold.
 Fortunately, $f_{*}(y)$ is strictly negative for $y<0$ and
 attracts no orbits from the positive quarter-plane.

Secondly, the necessary behaviour at infinity of $f(y)$ is:
 \be
 \label{in1}
  \tex{
 f(y)=C_0 |y|^{\frac {3\a}{1+\a}}(1+o(1)) \asA y \to - \iy \quad
 \big( \frac{3\a}{1+\a}=\frac \a \b \big),
 }
  \ee
  where $C_0>0$ is a constant, which can be arbitrarily changed by
  scaling \ef{2.8}.
 It is important for the future conclusions to derive the whole 3D
 bundle of solutions satisfying \ef{in1}. As usual, this is done
 by the linearization as $y \to -\iy$:
  \be
  \label{in12}
   \begin{split}
   & f(y)=f_0(y)+Y(y) \whereA f_0(y)=C_0(-y)^{\frac \a \b} \ssk\\
   \LongA
    &
     \tex{
      C_0 [(-y)^{\frac \a \b}Y]'''+ \b Y'(-y)+ \a Y + \frac 12\,
   [f_0(y)]'''+...=0.
   }
    \end{split}
   \ee
   According to classic WKBJ-type asymptotic techniques in ODE
   theory, we look for solutions of \ef{in12} in the exponential
   form with the following characteristic equation:
 \be
 \label{in22}
  \tex{
  Y(y) \sim {\mathrm e}^{a(-y)^\g}, \,\,\, \g=1+ \frac 12\,\big(1-
  \frac \a \b\big)>1
 \LongA  C_0(\g a)^2=-\b, \,\,
  a_\pm ={\rm i}\, \frac 1\g \, \sqrt{\frac {\b} C_0}.
  }
  \ee
This gives the whole 3D bundle of the orbits \ef{in1}: as $y \to
-\iy$,
 \be
 \label{in23}
  \tex{
  f(y)=C_0(-y)^{\frac \a \b}+
(-y)^\d \big[ C_1 \sin \big( \frac 1\g \, \sqrt{\frac {\b}
C_0}(-y)^\g \big)+  C_2 \cos \big( \frac 1\g \, \sqrt{\frac {\b}
C_0}(-y)^\g \big)\big]+...\, ,
 }
 \ee
where $ C_{1,2} \in \re$. The slow decaying factor $(-y)^\d$ in
 the double scale asymptotics \ef{in23} is not essential in what
 follows, so we do not specify the exponent $\d<0$ therein.

   The behaviour \ef{in22} and \ef{in23} give  crucial for us asymptotics: due to
  \ef{in1}, we have the gradient blow-up behaviour at a single
  point: for any fixed $x<0$, as $t \to 0^-$, where $y=x/(-t)^\b
  \to - \iy$,
   \be
   \label{in2}
    \tex{
    u_-(x,t)= (-t)^\a f(y) = (-t)^\a C_0 \big| \frac
    x{(-t)^\b}\big|^{\frac \a \b}(1+o(1)) \to C_0
    |x|^{\frac{3\a}{1+\a}},
    }
    \ee
    uniformly on compact subsets,
    as required by \ef{s3}.

 The geometry of the above asymptotic bundles yields:

\begin{proposition}
\label{Pr.1} For any fixed $A<0$ in $\ef{2.7}$, the problem
$(\ref{s2})$ admits the unique shock wave profile $f(y)$, which is
an odd  function and is strictly positive for $y<0$.
 \end{proposition}

Uniqueness follows from the asymptotics (\ref{2.7}) and scaling
invariance (\ref{2.8}).
 Global existence as infinite extension of
the unique solution from $y=0^-$ follows from the structure and
dimension of the bundles \ef{in23} (it is 3D, i.e., comprises the
whole phase space of the equation in this quadrant). The 2D ``bad"
bundle, composed from the orbits \ef{kk100}, is not connected with
the anti-symmetric ones.
 Concerning
 positivity, which is rather technical, see some details in \cite[\S~3,4]{GPnde}.

\begin{figure}
\centering
\includegraphics[scale=0.65]{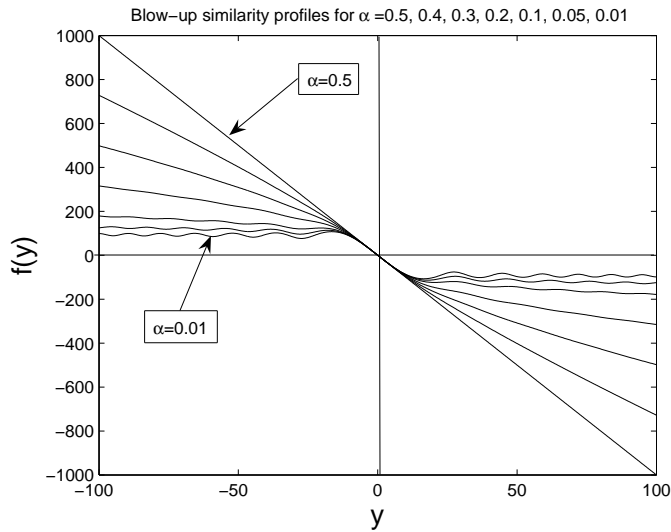}
\caption{\small The odd blow-up similarity profiles $f(y)$ in
$\re$ with $f'(0)=-10$ and $\a=0.5$, 0.4, 0.3, 0.2, 0.1, and
0.01.}
\label{F2}
\end{figure}
\begin{figure}
\centering
\includegraphics[scale=0.65]{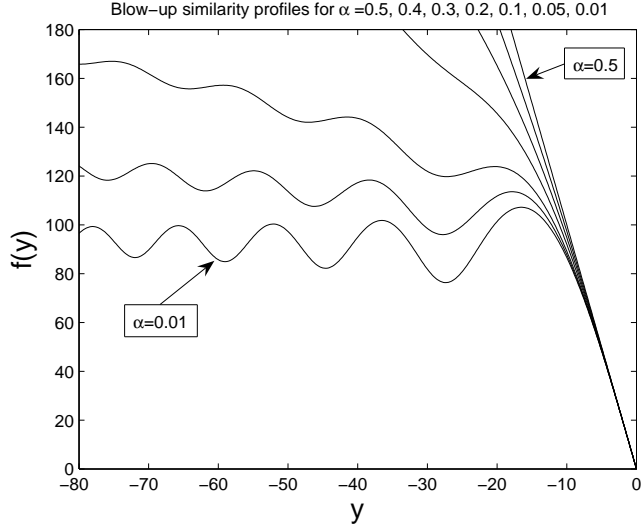}
 \vskip -.2cm
\caption{\small Enlarged oscillatory behaviour of $f(y)$ from
Figure \ref{F2} for $y \in [-80,0]$.}
   \vskip -.2cm
\label{F21}
\end{figure}
\begin{figure}
\centering
\subfigure[$y \sim -500$]{
\includegraphics[scale=0.52]{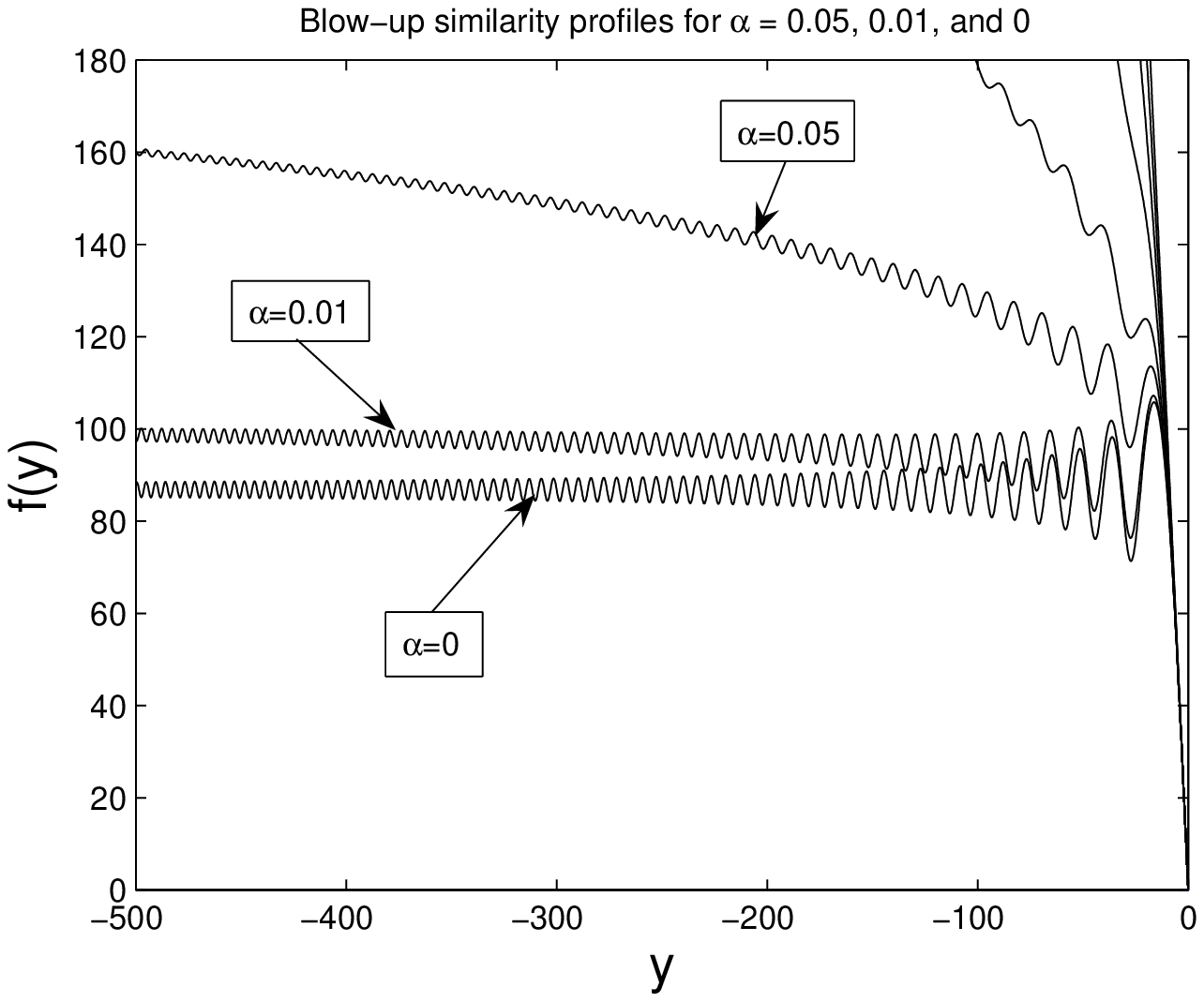}
}
\subfigure[$y \sim -1000$]{
\includegraphics[scale=0.52]{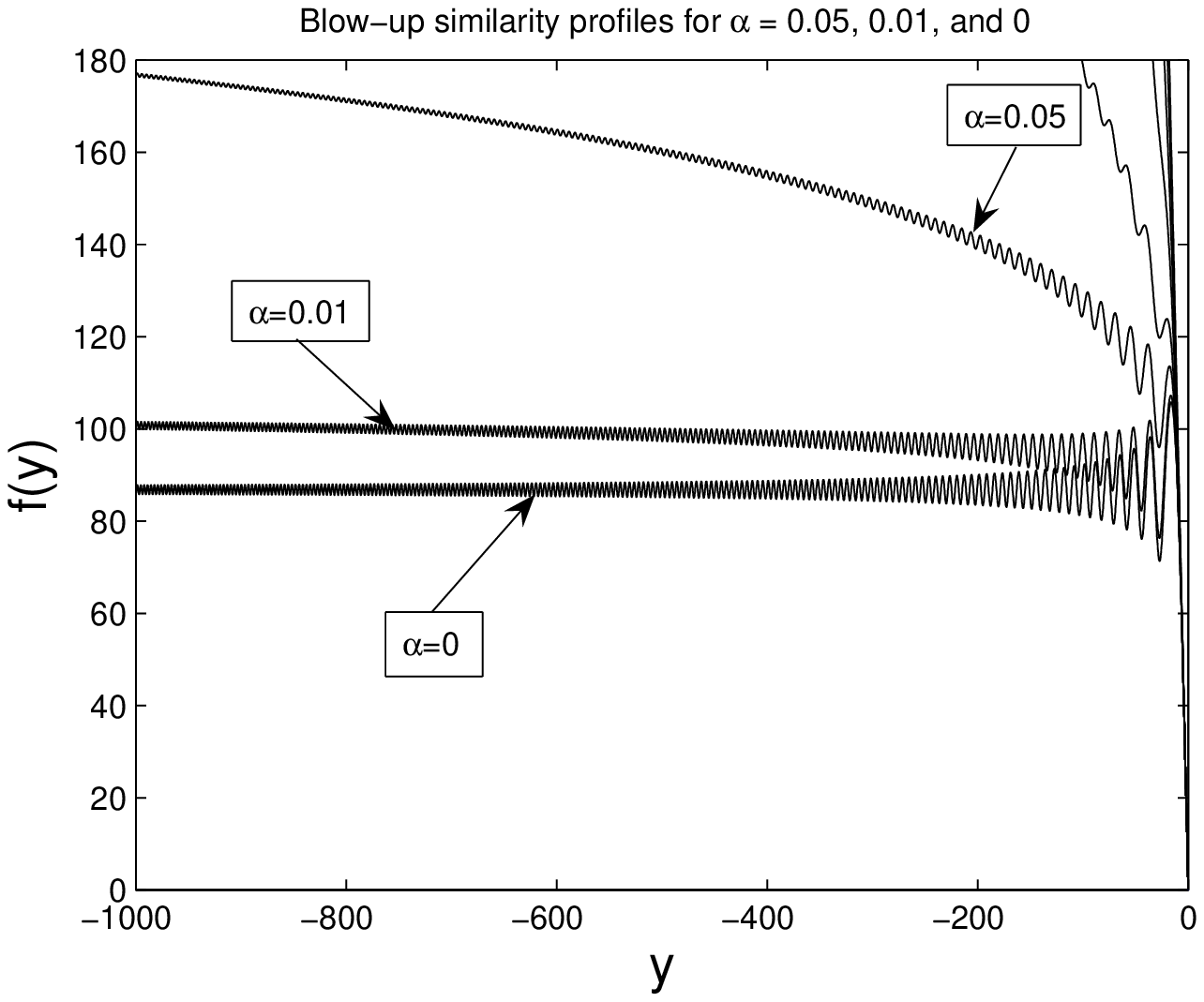}
}
 \vskip -.2cm
\caption{\rm\small Enlarged oscillatory behaviour of $f(y)$ from
Figure \ref{F2} for $y \in [-500,0]$ (a) and $y \in [-1000,0]$
(b).}
 \label{F22}
\end{figure}

 Figure \ref{F2} shows a general view of  similarity
profiles $f(y)$ for various values of the parameter $\a \in
[0.01,0.5]$ with the fixed derivative
 $$
 f'(0)=-10.
 $$
 In Figure \ref{F21}, we show enlarged oscillatory structure corresponding to \ef{in23} of
 some of the profiles closer to the origin. Figure \ref{F22}
 continues to explain the oscillatory behaviour \ef{in23} of solutions on
 different $y$-scales.

The next Figure \ref{F23} shows how the oscillatory features
dramatically increase for negative $\a$. For $\a = -0.099$, which
is very close to the critical value $- \frac 1{10}$ in \ef{ra1},
we observe a ``saw-type" profile of maximally allowed oscillatory
 structure; see more details in \cite[\S~4]{GPnde}.

\begin{figure}
\centering
\includegraphics[scale=0.65]{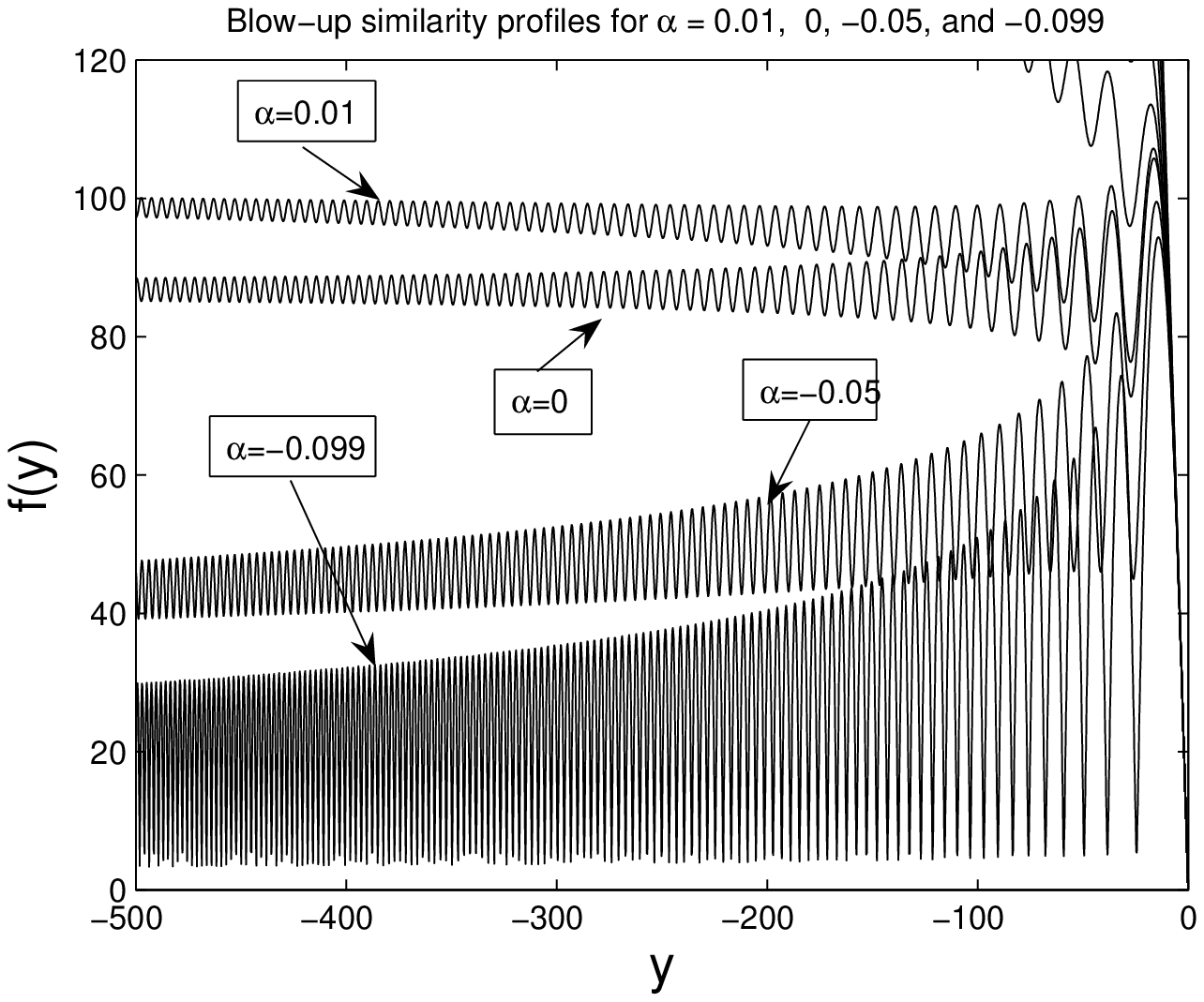}
 \vskip -.2cm
\caption{\small Oscillatory behaviour of $f(y)$ for positive and
negative $\a$.}
   \vskip -.2cm
\label{F23}
\end{figure}

Finally, we claim that, besides odd blow-up similarity profiles,
there exist other not that symmetric, for which the conditions at
the origin in \ef{s2} do not apply. This construction can be
performed similar to that in \cite[\S~3,4]{GPnde}. The blow-up
solutions constructed above are sufficient for our main purposes.

\subsection{On self-similar collapse of shocks}

It is curious that the similarity solutions \ef{s1} given by the
ODE \ef{s2} can describe {\em collapse} as $t \to 0^-$ of shocks.
These are given by profiles $f(y)$, which instead of the
anti-symmetry conditions in \ef{s2}, satisfy
 \be
 \label{col1}
 f(0)=f_0>0, \quad f'(0)=f_1<0, \andA f''(0)=f_2 \in \re.
  \ee
Since the bundle at infinity \ef{in23} is 3D and hence exhausts
all the trajectories there, shooting with the parameters \ef{col1}
yields an orbit $f(y;f_0,f_1,f_2)$, which, for  a wide range of
 $f_{0,1,2}$, has the behaviour \ef{in1} with a $C_0>0$ and
creates the data \ef{in2} as $t \to 0^-$.

In addition, for the corresponding similarity solution \ef{s1},
the shock disappears since
 \be
 \label{col2}
 [u_-(0,t)]= 2f_0(-t)^\a \to 0 \asA t \to 0^-.
  \ee
However, since the data at $t=0^-$ has a typical form \ef{in2}, a
new shock will be created as $t \to 0^+$, which we are going to
explain in Section \ref{S3}, where the parameterization such as in
\ef{col1} for the ODE \ef{s2N} will be used in greater detail.

\subsection{On smooth rarefaction waves}

These are global in time solutions of the form \ef{s1N}, with the
ODE as in \ef{s2N}. Then taking $F(y) \equiv -f(y)$, one observes
how the weakly singular initial data
 \be
 \label{in40}
 - u_-(x,0^-) \quad \mbox{given in (\ref{s3})}
 \ee
collapse into the smooth (even analytic) solution $u_+(x,t)$ for
$t>0$. This is quite similar to the same phenomena as in
\cite{GPnde} and \cite{GPndeII}, and we will not comment on this
anymore, but indeed will use the global solutions \ef{s1N},
however in their discontinuous hypostasis.

\section{Nonunique similarity extensions beyond blow-up}
 \label{S3}

\subsection{Nonuniqueness of similarity solutions}

As we mentioned, a discontinuous shock wave extension of blow-up
solutions \ef{s1}, \ef{s2} are performed by using the global ones
\ef{s1N}, \ef{s2N}. Actually, this leads to watching a whole
three-parametric family of solutions parameterized by their Cauchy
values at the origin:
 \be
 \label{N1}
 F(0)=F_0>0, \quad F'(0)=F_1<0, \andA F''(0)=F_2 \in \re.
  \ee

 The 3D phase space for the ODE in \ef{s2N} has two clear stable
 ``bad" bundles:

\ssk

 \noi{\bf (I)} positive solutions with ``singular extinction" in finite $y$, where
 $F(y) \to 0$ as $y \to y_0^+<0$. Indeed, this is an unavoidable
 singularity following from the degeneracy of the equations with
 the principal term $F F'''$ leading to the singular potential $\sim \frac
 1F$;

\ssk

 \noi{\bf (II)} positive solutions with the fast growth about the explicit solution:
  \be
  \label{f1}
   \tex{
  F_*(y) = - \frac {y^3}{60} \to +\iy \asA y \to - \iy.
  }
  \ee
 The 2D stable bundle is similar to that in \ef{kk100}.

\ssk

 Both sets of such solutions are {\em open} by the standard
continuous dependence of solutions of ODEs on parameters. The
desired solutions are situated in between those two stable open
bundles; cf. arguments in \cite[\S~3,4]{GPnde}.

As usual, it is key  to derive the whole 2D
 bundle of solutions satisfying \ef{in1}.
 Similar to \ef{in12}, we perform the standard  linearization as $y \to -\iy$ in \ef{s2N}:
  \be
  \label{in12G}
   \begin{split}
   & f(y)=F_0(y)+Y(y) \whereA F_0(y)=C_0(-y)^{\frac \a \b} \ssk\\
   \LongA
    &
     \tex{
      C_0 [(-y)^{\frac \a \b}Y]'''- \b Y'(-y)- \a Y + \frac 12\,
   [F_0(y)]'''+...=0.
   }
    \end{split}
   \ee
   The WKBJ method now leads to a different characteristic equation:
 \be
 \label{in22G}
  \tex{
  Y(y) \sim {\mathrm e}^{a(-y)^\g}, \,\,\, \g=1+ \frac 12\,\big(1-
  \frac \a \b\big)>1
 \LongA C_0(\g a)^2=\b, \,\,
  a_\pm = \pm \frac 1\g \, \sqrt{\frac {\b} C_0},
  }
  \ee
 so that the only admissible root is $a_-<0$.
This gives a 2D bundle of the orbits \ef{in1}:
 \be
 \label{in23G}
  \tex{
  f(y)=C_0(-y)^{\frac \a \b}+
(-y)^\d C_2 {\mathrm e}^{a_-(-y)^\g} +... \asA y \to -\iy, \quad
C_1 \in \re.
 }
 \ee

As a result, we have the following:

\begin{proposition}
\label{Pr.2}
 Let $\a \in (0, \frac 12)$.
 For any fixed $F_0>0$ and $F_1<0$ in $\ef{N1}$,
 there exists a unique $F_2 \in \re$ such that
  the problem
$(\ref{s2N})$ has a solution $F_*(y)$ for some $C_0>0$.
 \end{proposition}

The proof is performed by shooting as in \cite{GPnde, GalNDE5} by
using the stable bundles indicated in (i) and (ii) above. Such
techniques are currently well-established for various higher-order
ODEs. As a similar and more complicated example of a fourth-order
ODE, we refer to the  methods in \cite{Gaz06}, where by shooting
technique existence and uniqueness of a positive solution of the
radial bi-harmonic equations with source:
 \be
 \label{rrr1}
 \D^2_r u= u^p \forA r=|x|>0, \quad u(0)=1, \quad u'(0)=u'''(0)=0, \quad
 u(\iy)=0,
  \ee
  was proved in the supercritical Sobolev range
   $$
   \tex{
   p > p_{\rm Sob}= \frac{N+4}{N-4} \whereA N>4.
   }
   $$
   Here, analogously, there exists a single shooting parameter
   being the second derivative at the origin $F_2=u''(0)$; the value $F_0=u(0)=1$ is
    fixed by obvious scaling.
 Proving uniqueness of such a solution in \cite{Gaz06} is not easy and lead
 to essential technicalities, which the attentive Reader can
 consult in case of necessity.
 Note that, instead of the global behaviour such as \ef{f1}, the equation
 \ef{rrr1} admits the blow-up one governed by the principal
 operator
  $$
  u^{(4)}+...= u^p \quad (u \to +\iy).
   $$
The solutions vanishing at finite point otherwise can be treated
as in the family {\bf (I)}.


\ssk

\begin{figure}
\centering
\includegraphics[scale=0.45]{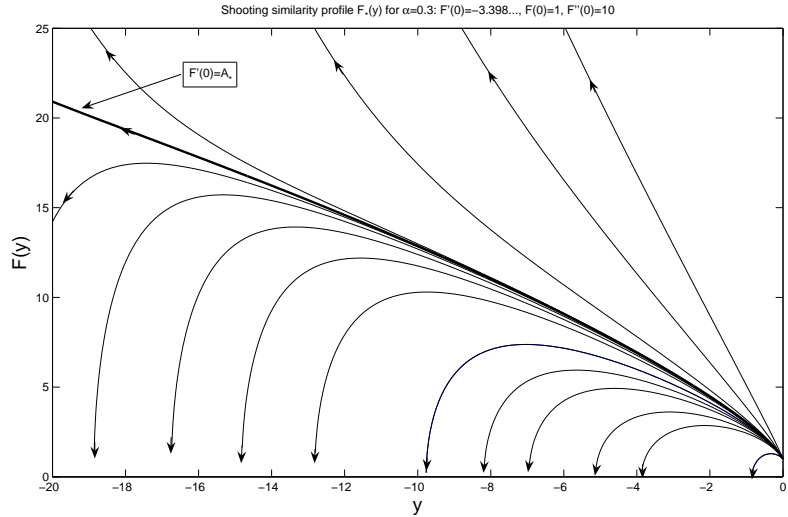}
 \vskip -.2cm
\caption{\small Shooting a proper solution  $F_*(y)$ of \ef{s2N}
for $\a=0.3$ and data $F(0)=F_0=1$, $F''(0)=F_2=10$, with
$F_1=F'(0)=A_*=-3.398...$.}
   \vskip -.2cm
\label{F5}
\end{figure}

\begin{figure}
\centering
\includegraphics[scale=0.65]{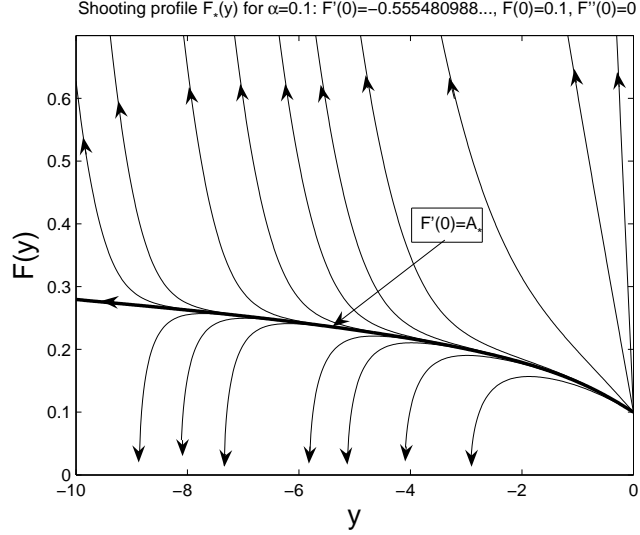}
 \vskip -.2cm
\caption{\small Shooting a proper solution  $F_*(y)$ of \ef{s2N}
for $\a=0.1$ and data $F(0)=F_0=0.1$, $F''(0)=F_2=0$, with
$F_1=F'(0)=A_*=-0.55548098...$.}
   \vskip -.2cm
\label{F6}
\end{figure}

\begin{figure}
\centering
\subfigure[$F_0=1$, $F_1=-1$]{
\includegraphics[scale=0.52]{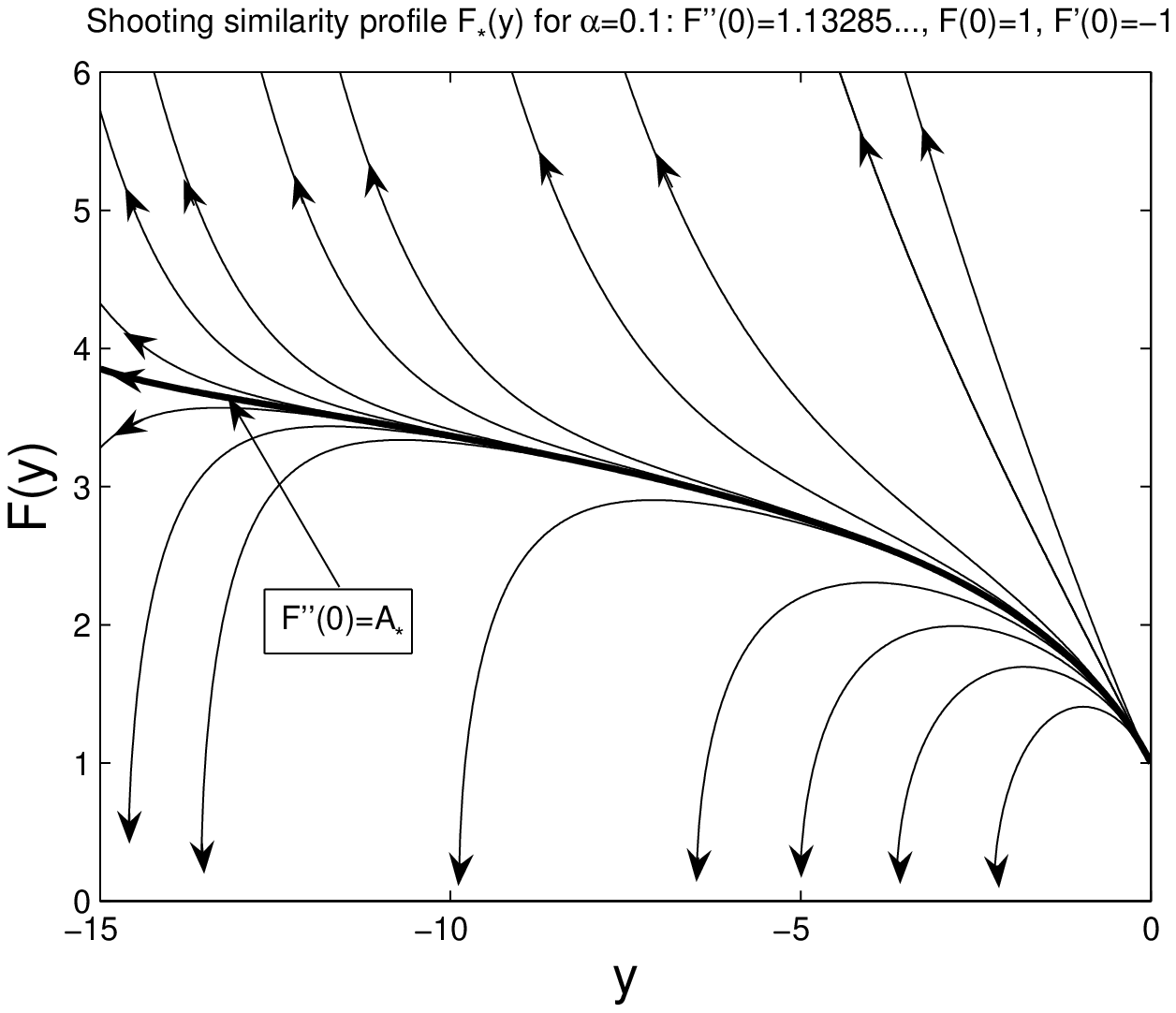} 
}
\subfigure[$F_0=1$, $F_2=0$]{
\includegraphics[scale=0.52]{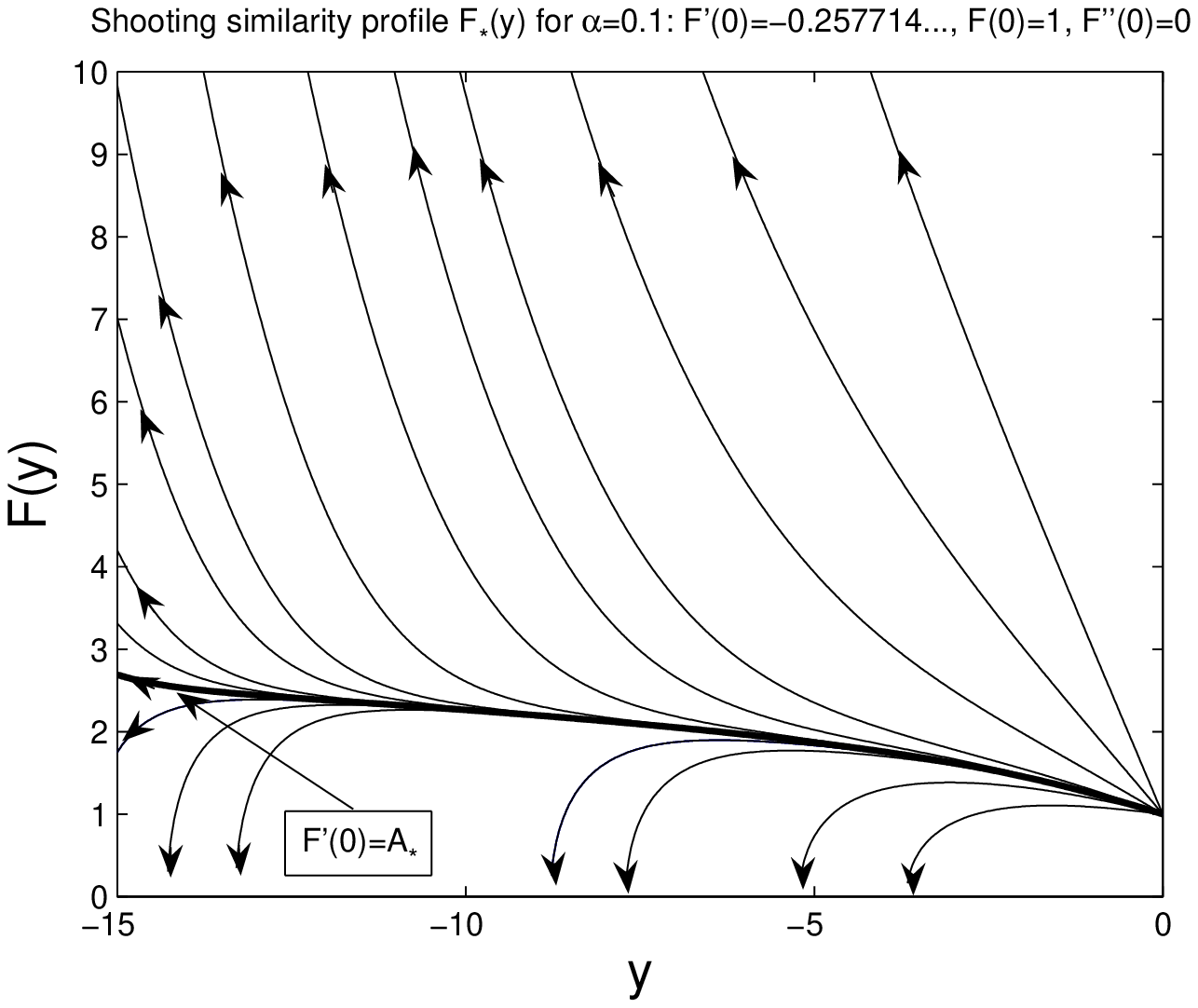} 
}
 \vskip -.2cm
\caption{\rm\small Shooting a proper solution  $F_*(y)$ of
\ef{s2N} for $\a=0.1$ and data $F(0)=F_0=1$, $F'(0)=F_1=-1$,
$F''(0)=F_2=A_*=1.13285...$ (a) and  $F(0)=F_0=1$, $F''(0)=F_2=0$,
$F'(0)=F_1=A_*=-0.257714...$  (b).}
 \vskip -.2cm
 \label{F7}
\end{figure}

Thus, as a result, we obtain a two-parametric family of solutions
of \ef{s2N} with an arbitrary fixed $C_0>0$ parameterized by
$F_0>0$ and $F_1<0$. For a fixed constant $C_0>0$ (uniquely given
by the blowing up limit $t \to 0^-$), the family is
one-parametric:
 \be
 \label{rr2}
 \Phi_{C_0}=\{F_*(y;F_0), \,\, F_0>0\},
  \ee
  which for convenience we parameterize by the value $F_0=F_*(0)$
  that measures the jump of the shock.
 The actual parameterization of the family \ef{rr2} is not of
importance. For instance, in Figure \ref{F5}, we show shooting
$F_*(y)$ for $\a=0.3$, which actually explains the strategy of
proving of Proposition \ref{Pr.2}. In Figure \ref{F6}, the same is
done for $\a=0.1$ and smaller $F(0)=0.1$, $F''(0)=0$. In Figure
\ref{F7}, again for $\a=0.1$, we show a couple of more shooting
with respect to parameters $F_2$ (a) and $F_1$ in (b).

\subsection{Remark: on nonexistence of solutions with $F_0=0$}

\begin{figure}
\centering
\subfigure[$\a=0.3$]{
\includegraphics[scale=0.52]{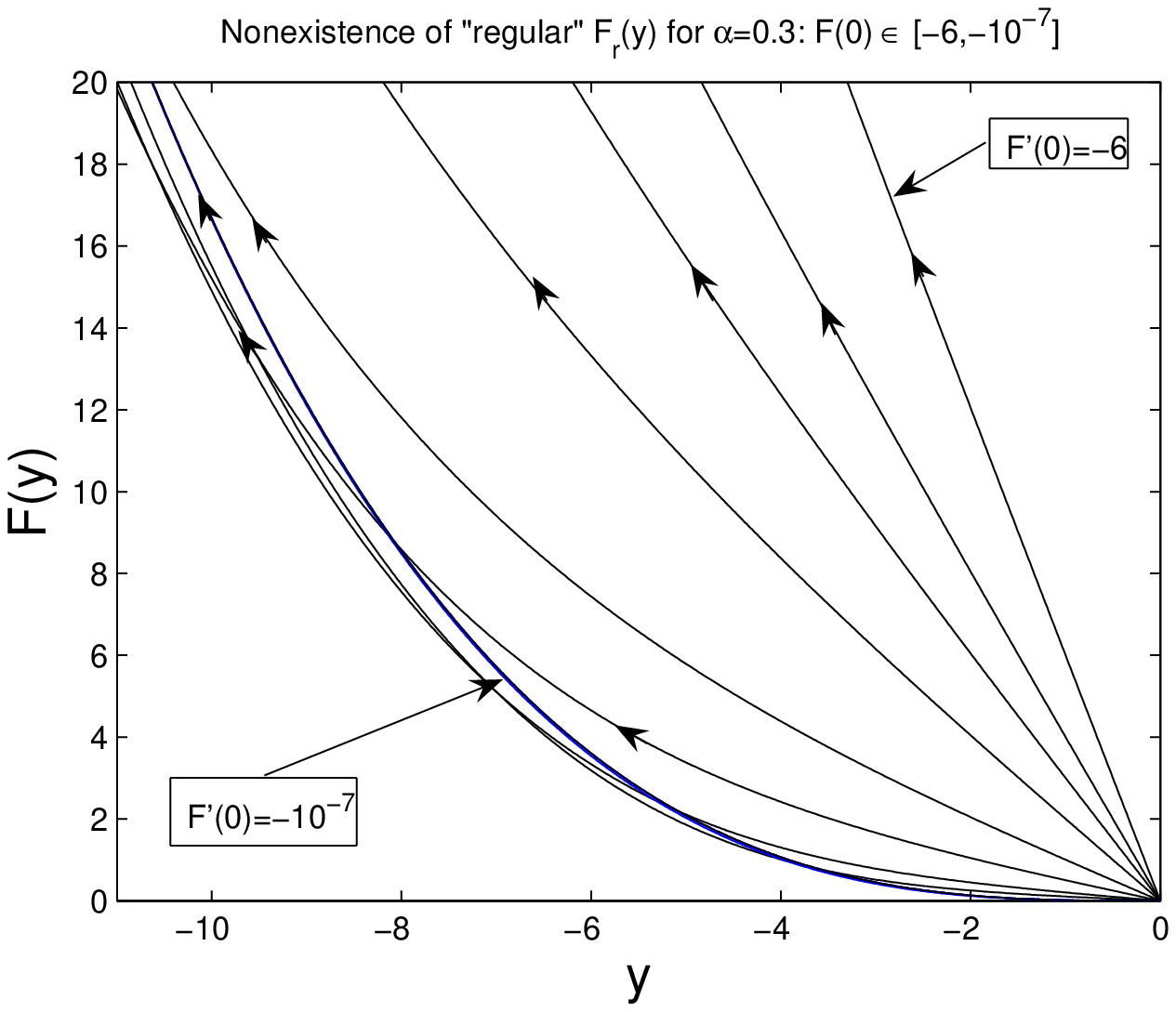} 
}
\subfigure[$\a=0.1$]{
\includegraphics[scale=0.52]{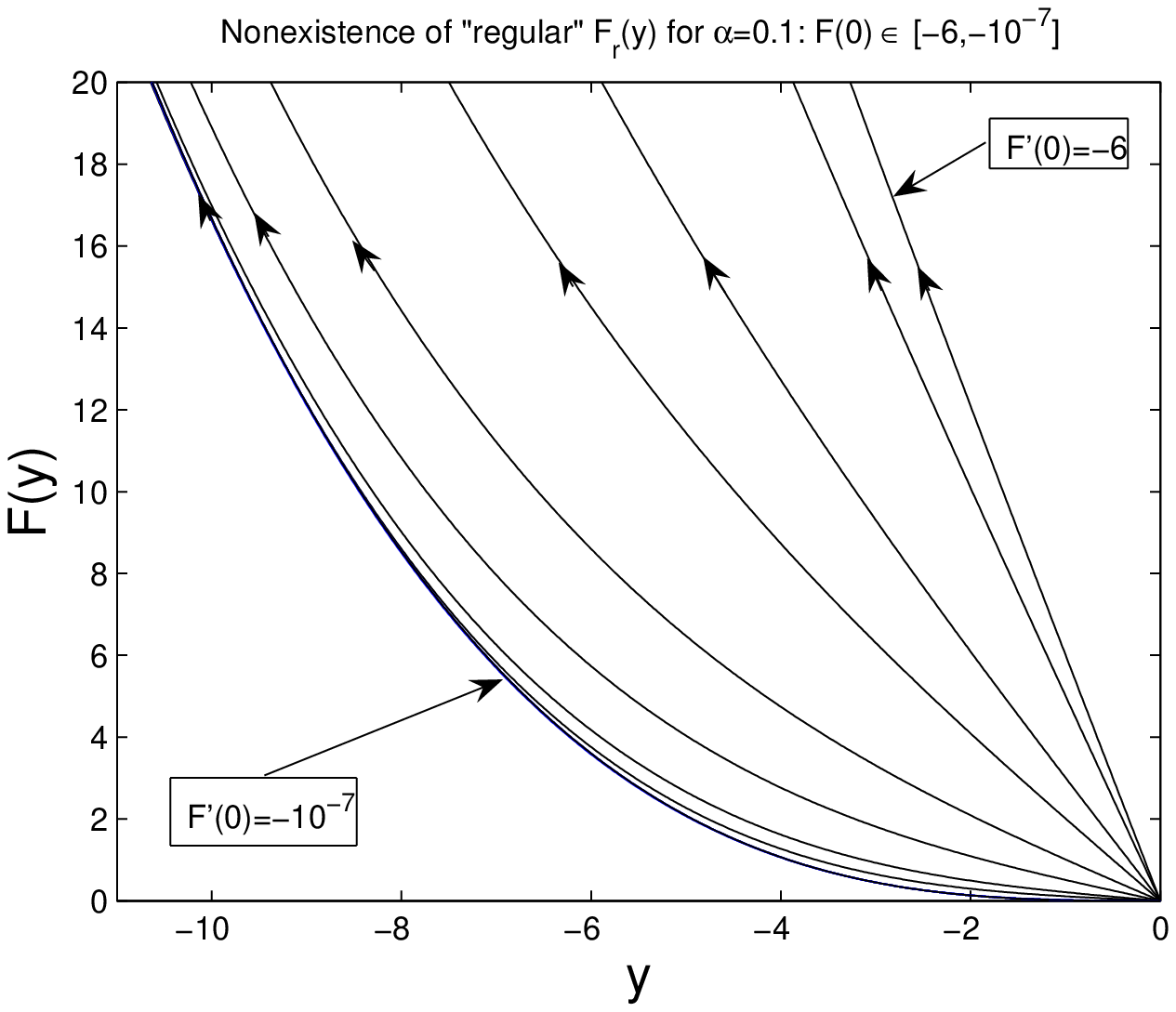} 
}
 \vskip -.2cm
\caption{\rm\small Towards  nonexistence of a solution of
\ef{s2N}, \ef{2.7N}, with $F_0=0$, for $\a=0.3$ (a) and $\a=0.1$
(b).}
 \vskip -.2cm
 \label{F8}
\end{figure}

This is also a principal issue. Indeed, if a solution of \ef{s2N}
for $F_0=0$, i.e., $F(y)$ having the asymptotics near the origin
similar to \ef{2.7},
\be
 \label{2.7N}
  \mbox{$
 F(y) = A y- \frac {1-2\a}{72}\, y^3+
  \frac {(1-2\a)^2}{72^2}\, \frac 1A \, y^5 +...\, \asA y \to 0
  \quad (A<0),
  $}
  \ee
 would exist, then the similarity solution \ef{s1N} would
 describe a smooth collapse of the initial singularity \ef{s3},
 and would actually mean an extra nonuniqueness in the problem.
 Fortunately, this is not the case and the single parameter $A<0$
 in \ef{2.7N} (actually reducing to $A=-1$ by scaling \ef{2.8})
 is not sufficient to shoot the necessary asymptotics as $y \to
 -\iy$ given in the second line in \ef{s1N}, i.e., the
 corresponding asymptotic bundles as $y \to 0^-$ and $y \to -\iy$
 are not overlapping. We do not prove this carefully (but indeed
 this can be done, since the ``non-overlapping of bundles" is
 large enough), and, as a key illustration, present Figure \ref{F8},
 where this nonexistence is carefully (with Tolerances about $10^{-10}$)
  checked for $\a=0.1$ (a) and $\a=0.3$ (b).


However, for more complicated nonlinear PDEs, such as, e.g.,
\cite{GalNDE5}
 $$
 u_{tt}=-(u u_x)_{xxxx}, \quad \mbox{etc.},
 $$
where the resulting ODEs are higher order with multi-dimensional
phase space, by some reasons, such solutions can be available,
meaning that formation of shocks can be reversible and is not a
unique option.


\subsection{Further discussion around nonuniqueness and non-entropy issues}

Thus, we have explained the required non-uniqueness \ef{s5N} of
the solution of the Cauchy problem \ef{1}, \ef{s3}: taking any
profile $F_* \in \Phi_{C_0}$ yields the self-similar continuation
\ef{s1N}, with the following behaviour of the jump at $x=0$:
 \be
 \label{jj1}
  \tex{
 -[u_+(x,t)]\big|_{x=0}
 \equiv - \big(u_+(0^+,t)-u_+(0^-,t))
 = 2F_0 t^\a>0 \forA t>0.
 }
 \ee
Rephrasing the result, let us emphasize that, in the similarity
(i.e., an ODE) representation, it has a pure dimensional origin:
the problem \ef{s2N} has too many solutions. More precisely, for
the given 3D shooting via parameters $\{F_0,F_1,F_2\}$ in the 3D
phase space of the ODE in \ef{s2N}, the dimension 1D of the
non-suitable
 bundles in {\bf(I)} and {\bf (II)} is not sufficient to define a
 unique solution $F_*(y)$, up to scaling \ef{2.8} as usual. In
 other words, the desired uniqueness could be achieved if
  \be
  \label{un55}
  \mbox{the ``bad" bundles in {\bf(I)} and {\bf (II)} are 2D in the 3D phase space}
   \ee
   (plus certain natural structural ``transversality" as a genetic property).
Hence, \ef{un55} is a geometric recognition of a possible
uniqueness (and entropy that assumes extra hard work) extension of
singular shock wave solutions.

Note that since these shocks are stationary, the corresponding
Rankine--Hugoniot (the R--H) condition on the speed $\l$ of the
shock propagation:
 \be
 \label{jj2}
 \tex{
 \l=- \frac {[(uu_x)_x]}{[u]}\big|_{x=0}=0
 }
  \ee
 is valid by anti-symmetry. As usual, this condition is obtained
 by integration of the equation \ef{1} in a small neighbourhood of
 the shock. Alternatively, \ef{jj2} is obtained by approximating
 the solution via a travelling wave (TW)
  \be
  \label{jj3}
   \tex{
  u(x,t)=f(x- \l t) \,\Longrightarrow \, -\l f'=(f f')'' \,\Longrightarrow \,
  -\l f=(f f')' \,\Longrightarrow \,  \l=- \frac {[(f f')']}{[f]}\big|_{y=0},
   }
   \ee
 which coincides with \ef{jj2}. Recall again that the R--H condition does
 not assume any novelty and is a corollary of integrating  the PDE
about the line of discontinuity.

Moreover, the R--H condition \ef{jj2} also indicates another
origin of nonuniqueness: a {\em symmetry breaking}. The point is
that the solution for $t>0$ is not obliged to be an odd function
of $x$, so we can define the self similar solution \ef{s1N} for
$x<0$ and $x>0$ using different triples of parameters
$\{F_0^\pm,F_1^\pm,F_2^\pm\}$, and the only extra condition one
needs is the R--H one:
 \be
 \label{kk1}
 \l=0 \LongA [(F F')']=0, \quad \mbox{i.e.,}
 \quad F_0^- F_2^- +(F_1^-)^2=   F_0^+ F_2^+ +(F_1^+)^2,
  \ee
 which can admit other solutions than the obvious anti-symmetric
 one
  $$
  F_0^-=-F_0^+, \quad F_1^-=F_1^+, \andA F_2^-=-F_2^+.
   $$

\ssk

 Overall, in any case, we have at least a one-parameter family of
solutions of the CP $\Psi_{C_0}=\{u_+(x,t)\}$ for $t>0$ depending
on arbitrary parameter  $F_0>0$. This family {\em is not discrete}
and hence does not reveal any particular solution $u_+$,
exhibiting some special properties as being, say, maximal,
minimal, or extremal in any sense. Under a natural assumption,
this confirms the extra negative statement \ef{s6N} on
nonexistence of a sufficiently general entropy-like inequality,
condition, and/or a procedure to detect a unique solution, as a
``better" special one. At least, if even an ``entropy-like"
procedure would have been derived (somehow, by any hypothetical
means), in view of the nonunique formation of shocks everywhere,
the resulting solution  would have never been  relative to data in
any metric, i.e., Hadamard's well-posedness concept would have
been violated anyway.

Indeed, different regularization of the NDE can lead to different
solutions, i.e., by using a parabolic regularization:
 \be
 \label{par1}
 u_\e: \quad u_t=(u u_x)_{xx}- \e u_{xxxx} \quad (\e>0),
  \ee
  with the same data $u_0(x)$.
Proving that the regularized sequence $\{u_\e, \e>0\}$ is a
compact subset in some metrics is indeed a very difficult problem,
which remains open for general data (see \cite{GalEng} for some
details concerning such problems). Anyway and however, assuming
that a suitable compactness of $\{u_\e\}$ is already available, we
do not think that a {\em proper unique solution} $\bar u(x,t)$ can
be obtained as the unconditional limit
 \be
 \label{ll1}
  \tex{
 \bar u(x,t)= \lim_{\e \to 0^+} \, u_\e(x,t).
  }
  \ee
More precisely, in view of the above nonuniqueness of post-blow-up
similarity extensions, we believe that \ef{ll1} has infinitely
many partial limits and, along various subsequences $\{\e_k\} \to
0^+$, the corresponding sequences $\{u_{\e_k}\}$ can converge to
different solutions $u_+$ with various values of rescaled shocks
$F_0>0$.
 We justify this in an example below.

 \ssk

 \noi{\bf Example.} We take $\a= \frac 15$ and $C_0=1$, so that
 the blow-up initial data \ef{s3} are
  \be
  \label{u1}
  u_0(x)= \sqrt {|x|} \forA x\le 0,
   \ee
   with the odd extension for $x>0$. Performing in \ef{par1} a natural scaling
 \be
 \label{u2}
  \tex{
  u_\e(x,t)= \e^{\frac 13} v(y,\t), \quad y= \frac x{\e^{2/3}},
  \quad \t= \frac t{\e^{5/3}}
  }
   \ee
   deletes $\e$ and
 yields the following uniformly parabolic problem:
  \be
  \label{u3}
  v_\t=(v v_y)_{yy} - v_{yyyy} \inB \re_- \times \re_+, \quad v_0(y)=
  \sqrt{|y|}\inB \re_-,
   \ee
   with the anti-symmetry conditions $v=v_{yy}=0$ at $y=0$.
 It is not a great deal, by using classic parabolic theory,
  to prove global existence of a unique
 classical solution of the CP \ef{u3}, though it takes some
 efforts. But the main issue is not existence/uniqueness. Indeed,
according to \ef{u2}, the study of the convergence \ef{ll1} as $\e
\to 0^+$ reduces to delicate asymptotic behaviour of $v(y,\t)$
simultaneously as $\t \to +\iy$ and $y \to -\iy$, which represents
a very difficult problem. However,  the above detected
nonuniqueness of the similarity extensions makes such an open
asymptotic problem excessive and not that necessary (quite
fortunately it seems). Note that studying the behaviour of the
hypothetical proper limit \ef{ll1} at the point of discontinuity
$x=0$ assumes detecting the limit of $\{u_\e(x_k,t)\}$ ($t>0$
small is fixed) for rather arbitrary sequences $\{x_k\} \to 0^-$
as $\e \to 0$. We then claim that such a limit essentially depends
of the choice of $\{\e_k\} \to 0^+$, i.e., the limit of the
sequence
 \be
 \label{u4}
 \tex{
 u_{\e_k}(x_k,t) \equiv \e_k^{\frac 13} v\big(\frac {x_k}{\e_k^{2/3}},
  \frac t{\e_k^{5/3}}\big) \asA k \to \iy
  }
  \ee
is very much $\{x_k\}$- and $\{\e_k\}$-dependent.

\ssk

We also expect that using other types of regularization, e.g., by
classic Bubnov--Galerkin methods of finite-dimensional
approximations on suitable functional Riesz  bases (see e.g.,
strong applications in Lions' classic monograph \cite{LIO}), which
obviously give globally existing solutions, will also lead to
similar nonuniqueness issues in the limit\footnote{Therefore,
 the issue of the {\em Galerkin uniqueness}, which can be treated
 as
 a simplest
   idea for
 an entropy-like construction,  becomes author's unrealizable dream.},
 so
these are unavoidable difficulties of modern PDE theory.

In other words, we then expect that, for a number of
higher-order\footnote{This is a principle issue: for general
second-order nonlinear parabolic equations with blow-up and
obeying the Maximum Principle, there always exists a unique {\em
proper minimal solution}, which does not depend on the type of
monotone regularization of the equation; see \cite[\S~2]{GV} or
\cite[Ch.~6]{GalGeom}.} nonlinear PDEs with singularities such as
blow-up, extinction, or shock wave formation, fully consistent
(see above) general entropy-like procedures to reveal a unique
solution cannot be available in principle. As we have shown, this
is just prohibited by a sufficient dimension of the phase space
(depending on the spatial order $\ge 3$ of the PDE), which allows
to shoot  a continuous family of solutions beyond singularity,
whose family does not have any isolated and/or  boundary points.

\subsection{On uniqueness for the FBP setting}

Evidently, the uniqueness can be restored if an extra condition at
the shocks is assumed, which poses an FBP for \ef{1}. For
instance, following Figure \ref{F6} and \ref{F7}(b), this happens
if we fix
 \be
 \label{hh1}
 F_2=F''(0)=0 \LongA u_{xx}(0^\pm,t)=0.
  \ee
Indeed, the uniqueness (in the present self-similar setting) is
restored if the set \ef{rr2} contains a unique such profile, i.e.,
 \be
 \label{hh2}
  \Phi_{C_0} \cap \{F_*''(0)=0\}= \{\hat F_*(y), \,\, \mbox{with a fixed}\,\,\hat F_0>0\}.
   \ee
Then this fixes the unique ``shock divergence" \ef{jj1} with
$F_0=\hat F_0$ beyond the singularity. In a PDE setting,
existence-uniqueness of a solution of the FBP \ef{1}, \ef{hh1} is
a difficult problem, which, for some simple geometric
configuration of shocks, can be solved by traditional FBP methods,
such as von Mises transformations and others. In general,
mathematical difficulties can be extremely challenging.

 A more general free-boundary condition on shock lines can be
 predicted from the structure \ef{s1N} of (generic) similarity
 continuation beyond blow-up.
 One should only take into account that parameters $\a$ and $\b$
 cannot enter such conditions, since these essentially depend on
 {\em a priori} unknown blow-up  ``initial data" $u(x,0^-)$.
  Since at the shock at $x=0$, by \ef{s1N}
   \be
   \label{hh3}
   u_+=t^\a F_0, \quad (u_+)_x=t^{\a-\b} F_1, \andA (u_+)_{xx}=t^{\a-2\b}
   F_2,
    \ee
    it is easy to reconstruct a general FBP condition at
    shocks, which suits for arbitrary $\a$:
    \be
    \label{hh4}
    u u_{xx} = \kappa (u_x)^2 \whereA \kappa \in \re.
     \ee
 Since both sides are of order $t^{2\a-2\b}$, this condition well
 corresponds to any of similarity formation of shocks (what
 happens for other configurations, is another delicate story).
Then \ef{hh1} is obtained for $\kappa=0$. Formally, $\kappa=\iy$
yields $F_1=0$, i.e., a kind of
 $$
  \mbox{a ``Neumann" FBP}: \quad u_x=0 \quad \mbox{at shocks}.
 $$

 Of course, it is necessary to check, for which $\kappa$ the
 condition
  \be
  \label{hh5}
  F_0 F_2= \kappa (F_1)^2
  \ee
  yields a unique profile $F_* \in \Phi_{C_0}$, and this occurs for arbitrary $\a$.
 As customary, posing necessary free-boundary conditions is an applied physical issue, though
  checking and predicting the well-posedness of the FBPs occurred is
  indeed a mathematical problem.

\section{Final remarks: on origin of uniqueness for the 1D Euler
equation}
 \label{S.4}

In this connection, it is interesting to trace out similar origins
of the uniqueness in the CP for the Euler equation \ef{3}. Its
obvious advantage is that it is solved via characteristics and the
solutions are given by an algebraic relation:
 \be
 \label{al1}
  \tex{
   \frac {{\mathrm d}t}1= \frac{{\mathrm d}x}u \LongA x-u \, t={\rm const.} \LongA
 u(x,t)= u_0(x- u(x,t)t).
 }
 \ee
 Let us assume the same initial data \ef{s3}, so that
 for $x \le 0$ we have the equation
  \be
  \label{al2}
   \tex{
   u=C_0(u \, t-x)^{\frac \a \b}, \quad \mbox{where now, dimensionally,} \quad \b=1+\a.
   }
   \ee
 Obviously, setting here $x=0$, corresponding to the permanent
 position of this stationary shock, yields the {\em unique} value
 of the solution at the shock:
\be
  \label{al3}
   \tex{
   u=C_0(u t)^{\frac \a \b} \LongA u(0^-,t)= C_0^{ \frac \b{\b-\a}}t^{\frac
   \a{\b-\a}} \equiv C_0^\b t^\a \forA \b=1+\a.
   }
   \ee
Or, analogously, the same uniqueness is guaranteed by the fact
that, for the first-order ODEs corresponding to rescaling of
\ef{3}, the phase space is also one-dimensional, that allows a
unique (and very simple) matching of two bundles. The
corresponding similarity solutions for $t<0$ and $t>0$
respectively are given by:
 \be
 \label{sim1}
  \begin{matrix}
  u_-(x,t)=(-t)^\a f(y), \,\,\, y= \frac x{(-t)^\b}
  \LongA f f'+ \b f'y -\a f=0,\qquad\qquad \ssk\ssk\\
\quad u_+(x,t)=t^\a F(y),\,\, \,\,\, y= \frac
x{t^\b}\quad\quad\,\,
  \LongA F F'- \b F'y +\a F=0,\qquad\qquad
  \end{matrix}
   \ee
   where $\b = 1+\a$,
and then as in  \ef{al3} we have to have
  \be
  \label{ha1}
 \tex{
u_-(0^-,t)= C_0^{\b}\, t^{\a} \forA t>0, \,\,\, \mbox{i.e.,}
\,\,\, F(0)=C_0^\b .
  }
  \ee

Of course, the ODEs in \ef{sim1} are explicitly integrated. As
usual, the last one for the post-blowing up behaviour, is
responsible for a unique continuation beyond blow-up. Integrating
it yields the unique solution $F_*(y)$ defined by the algebraic
equation
 \be
 \label{un1}
 \tex{
  F'= - \frac {\a F}{F-\b y}, \,\,\, F=-yP, \,\,\, yP'= - \frac{P(P+1)}{P+\b} \LongA
  \frac{ F^\b(y)}{[F(y)+|y|]^\a}= C_0^\b \forA y\le 0,
  }
  \ee
  where $C_0$ is the constant in \ef{al2} (a similar procedure applies to the blow-up profile $f$). It then follows from \ef{un1} that
  the rescaled shock jump is also uniquely determined:
 \be
 \label{un2}
  F_0=F_*(0)=C_0 ^\b,
  \ee
  and indeed this precisely gives the time behaviour \ef{al3} and \ef{ha1} of the
  shock value.

  Note that, unlike calculus in \ef{kk1}, a symmetry breaking is
  obviously impossible here. Indeed, the R--H condition uniquely
  implies the symmetry:
   \be
   \label{rh1}
    \tex{
    \l=\frac {[F^2]}{2[F]}
 = \frac{(F_0^+)^2-(F_0^-)^2}{2(F_0^+-F_0^-)}
     \equiv \frac {F_0^-+F_0^+}2=0 \LongA
    F_0^+=-F_0^-.
    }
    \ee

It is not an exaggeration to say that this
 unique micro-structural extension of generic blow-up
 singularities at any point
eventually reflects a true dimensional and ODE origin of those
classic entropy theories constructed for conservation laws in
$\re$ by Oleinik in the 1950s and by Kruzhkov in 1960s for
equations in $\ren$. Indeed, if, for the simplest ODE similarity
problem for \ef{3}, the matching would be nonunique and
non-discrete (plus something else), these very influential
theories would have never been appeared and  created.
 We thus
claim that this nonuniqueness and non-entropy  are precisely  the
case for the Cauchy problem for the NDE \ef{1} and can be expected
for a number of other nonlinear dispersion (and not only those)
higher-order
 equations from PDE theory of the twenty first
century.



\begin{thebibliography} {44}



\bibitem
 {Bres}
 A.~Bressan, {\rm Hyperbolic Systems of Conservation Laws. The One
 Dimensional Cauchy Problem}, Oxford Univ. Press, Oxford, 2000.





\bibitem
 {Chr07}
  D.~Christodoulou,{\em The Euler equations of compressible fluid
  flow}, Bull. Amer. Math. Soc.,
 {\bf 44} (2007), 581--602.



\bibitem
 {Daf}
 C.~Dafermos, {Hyperbolic Conservation Laws in Continuum Physics},
 Springer-Verlag, Berlin, 1999.



\bibitem 
 {GalGeom}
 V.A.~Galaktionov, {\rm Geometric Sturmian  Theory of Nonlinear
 Parabolic Equations and Applications}, Chapman$\,\&\,$Hall/CRC, Boca Raton,
Florida,
 2004.




\bibitem 
{GalEng} V.A.~Galaktionov, {\em  On higher-order viscosity
approximations
  of  odd-order nonlinear PDEs}, J.~Engr. Math.,
  {\bf 60} (2008), 173--208.


    \bibitem
{GPndeII}
 V.A.~Galaktionov, {\em
 Nonlinear dispersion equations: smooth deformations, compactons, and extensions to higher orders},
 {Comput. Math. Math. Phys.,} {\bf 48} (2008), 1823--1856
 (arXiv:0902.0275).



 \bibitem
 {GalNDE5}
 V.A.~Galaktionov,
 {\em Shock waves and compactons for fifth-order nonlinear dispersion
 equations}, Europ. J.~Appl. Math., submitted (arXiv:0902.1632).




  \bibitem
{GPnde}
 V.A.~Galaktionov and S.I.~Pohozaev, {\em
 Third-order nonlinear dispersive equations: shocks, rarefaction,
 and blow-up waves},
 {Comput. Math. Math. Phys.,} {\bf 48} (2008), 1784--1810
 (arXiv:0902.0253).




\bibitem 
{AMGV}
 V.A.~Galaktionov and J.L.~Vazquez, {A Stability Technique  for Evolution
Partial Differential Equations.
 A Dynamical Systems Approach},
Birkh\"auser, Boston/Berlin, 2004.

\bibitem 
{GV}
  V.A.~Galaktionov and J.L.~Vazquez, {\em Continuation of
blow-up solutions of nonlinear heat equations in several space
dimensions,} {Comm. Pure Appl. Math.,} {\bf 50} (1997), 1-68.





  \bibitem
{Gaz06}
 F.~Gazzola and H.-C.~Grunau,  {\em Radial entire solutions for
 supercritical biharmonic equations},
 Math. Ann.,
 {\bf 334} (2006), 905--936.




\bibitem 
 {Kru2}  {S.N.~Kruzhkov},
{\em First-order quasilinear equations in  several independent
variables}, {Math. USSR Sbornik}, {\bf 10} (1970), 217--243.




\bibitem 
{LIO}  J.L.~Lions, {\rm Quelques m\'{e}thodes de r\'{e}solution
des probl\`{e}mes aux limites non lin\'{e}aires\/}, Dunod,
Gauthier--Villars, Paris, 1969.




\bibitem
 {Ol1}
 O.A.~Oleinik, {\em Discontinuous solutions of non-linear
 differential equations}, Uspehi Mat. Nauk., {\bf 12} (1957), 3--73;
 Amer. Math. Soc. Transl. (2), {\bf 26}
 (1963), 95--172.

 \bibitem
  {Ol59}
  O.A.~Oleinik, {\em Uniqueness and stability of the generalized solution of the
  Cauchy problem for a quasi-linear equation}, Uspehi Mat. Nauk., {\bf 14} (1959), 165--170;
Amer. Math. Soc. Transl. (2), {\bf 33}
 (1963), 285--290.


  \bibitem
 {Poi08}
 D.~Poisson, {\em M\'emoire sur la th\'eorie du son}, J.~Polytech.
 (14 \'eme cahier) {\bf 7} (1808), 319--392.



\bibitem
 {Pom08}
  Y.~Pomeau, M.~Le~Berre, P.~Guyenne, and S.~Grilli,
  {\em Wave-breaking and generic singularities of nonlinear hyperbolic
  equations},
 Nonlinearity, {\bf 21}
  (2008), T61--T79.



\bibitem
 {Ri58}
  B.~Riemann, {\em \"Uber die Fortpfanzung ebener Luftwellen von
  endlicher Schwingungswete,  Abhandlungen der Gesellshaft der
  Wissenshaften zu G\"ottingen}, Meathematisch-physikalishe
  Klasse,
 {\bf 8} (1858-59), 43.








\bibitem  
{Sm}    J.~Smoller,
 {\rm Shock Waves and Reaction-Diffusion Equations,}
Springer-Verlag, New York, 1983.



\end{thebibliography}
\enddocument